\documentclass[11pt,letterpaper]{article}
\usepackage{ifthen,latexsym,amssymb,amsmath,amsthm}
\usepackage{mathrsfs}
\usepackage{graphicx}
\usepackage{tikz}
\usepackage{booktabs,blkarray,bigstrut}
\usepackage{url}
\usepackage{xskak,chessboard}

\usepackage{xcolor}
\definecolor{mylinkcolor}{rgb}{0.1, 0.5, 0.69} 
\usepackage[colorlinks=true,linkcolor={mylinkcolor}]{hyperref}
\hypersetup{
    citecolor = mylinkcolor
}

\newcommand{\ZB}[1]{\overline{B}_{#1}}
\newcommand{\coef}[2]{C^{#1}{#2}}

\newcommand{\QQ}{\mathcal{Q}}
\newcommand{\LL}{\mathcal{L}}

\newif\ifnotesw\noteswtrue
\noteswfalse   

\newcommand{\beq}[1]{\begin{equation}\label{eq:#1}}
\newcommand{\eeq}{\end{equation}}

\newtheorem{theorem}{Theorem}
\newcommand{\bth}[2][nothing]{\ifthenelse{\equal{#1}{nothing}}
 {\begin{theorem}} {\begin{theorem}[#1]}\label{th:#2}}

\newtheorem{lemma}[theorem]{Lemma}
\newcommand{\blm}[2][nothing]{\ifthenelse{\equal{#1}{nothing}}
 {\begin{lemma}} {\begin{lemma}[#1]}\label{lm:#2}}

\newtheorem{problem}[theorem]{Problem}
\newcommand{\bpr}[2][nothing]{\ifthenelse{\equal{#1}{nothing}}
 {\begin{problem}} {\begin{problem}[#1]}\label{pr:#2}}

\newtheorem{proposition}[theorem]{Proposition}

\renewcommand{\qed}{\nolinebreak\mbox{\hspace{5 true pt}%
  \rule[-0.85 true pt]{3.9 true pt}{8.1 true pt}}}

\newcommand{\nothree}{no-3-in-a-line }

\newcommand{\wangdata}{Xianzhi Wang \thanks{Department of Mathematics, Texas A\&M University, College Station, TX 77843-3368.}}

\newcommand{\ohdata}{Seunghwan Oh\thanks{Department of
    Mathematics, Middlebury College, Middlebury, VT 05753.}}

\newcommand{\schmittdata}{John R. Schmitt\thanks{Department of
    Mathematics, Middlebury College, Middlebury, VT 05753.  Supported in part by Fulbright Austria.  Corresponding author: jschmitt$@$middlebury.edu.}}

\newcommand\labelfive[1]{%
 \ifcase#1 \or -2\or -1\or 0\or 1\or 2\fi}

\newcommand\labelnine[1]{%
 \ifcase#1 \or -4\or -3\or -2\or -1\or 0\or 1\or 2\or 3\or 4\fi}

\title{Repeatedly applying the Combinatorial Nullstellensatz for Zero-sum Grids to Martin Gardner's minimum no-3-in-a-line problem}
\author{\ohdata, \schmittdata, \wangdata}
\date{}

\begin{document}

\maketitle
\begin{abstract}
  In 1976 Martin Gardner posed the following problem: ``What is
  the smallest number of [queens] you can put on an [$n \times n$ chessboard]
  such that no [queen] can be added without creating three in a row, a
  column, or a diagonal?''  The work of Cooper, Pikhurko, Schmitt and Warrington showed that this number is at least $n$, except in the case when $n$ is congruent to $3$ modulo $4$, in which case one less may suffice.    When $n>1$ is odd, Gardner conjectured the lower bound to be $n+1$.  We prove this conjecture in the case that $n$ is congruent to 1 modulo 4.  The proof relies heavily on a recent advancement to the Combinatorial Nullstellensatz for zero-sum grids due to Bogdan Nica.
\end{abstract}

\section{Introduction.}

The 1976 Mathematical Games column in {\it Scientific American} (see Gardner \cite{gardner:klondike}, Chapter 5,
pg. 71) posed the question: What is the minimum number of counters that can be placed on an $n\times n$ chessboard, no three in a line, such that adding one
more counter on any vacant square will produce three in a line?  This is the {\it minimum \nothree problem}.  There are different ways to interpret `line'.  We may take the word in its broadest sense, allowing for a straight line of any orientation.  In this case, we point the reader to the recent work of Aichholzer, Eppstein, and Hainzl \cite{aichholzereppsteinhainzl}.  However, we restrict our attention to when `line' refers only to orthogonals and diagonals (as was also done in \cite{cooperpikhurkoschmittwarrington:14}); this is the \emph{queens version} of the problem since it properly relates to how a queen is permitted to move in a game of chess.

\begin{figure}
    \centering
    \input{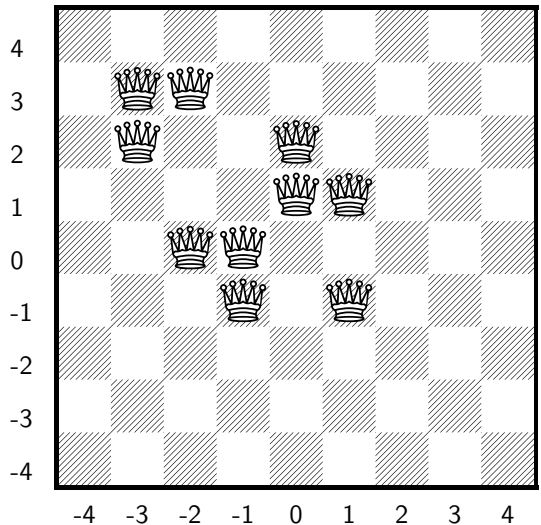}
    \caption{An example of a maximal placement with $10$ queens on a $9 \times9$ board.}
    \label{fig:10on9by9}
\end{figure}

Figure~\ref{fig:10on9by9} shows a $9\times 9$ chessboard with a placement of $10$ queens with no three in a line.  This placement is \emph{maximal}, that is, any additional queen will create three in a line.  For example, an additional queen placed in the top left corner produces three in a line along the main diagonal.  This particular placement is also of minimum size (where {\it size} of a placement is the number of queens in the placement), that is, there is no placement with nine or fewer queens meeting the requirements.


A history of the queens version problem is given in \cite{cooperpikhurkoschmittwarrington:14}\footnote{The published version of this paper contains an editorial error in the abstract.  The arxiv version 2 corrects this error.}.  Here we offer a brief summary of known results.  Cooper, Pikhurko, Schmitt and Warrington \cite{cooperpikhurkoschmittwarrington:14} gave the following.

\begin{theorem}\label{theorem:CPSW14}[Cooper, Pikhurko, Schmitt, Warrington \cite{cooperpikhurkoschmittwarrington:14}]
  For $n \geq 1$, the answer to Gardner's \nothree queens version problem is at least $n$, except in the case when $n$ is congruent to $3$ modulo $4$, in which case one less may suffice.
\end{theorem}

An elementary, {\it ad hoc} proof in the case of $n$ even was given, yet the approach only yielded a lower bound of $n-1$ when $n$ is odd.  To obtain the full result, the proof of Theorem \ref{theorem:CPSW14} relied on Alon's Combinatorial Nullstellensatz (reproduced here as Theorem~\ref{theorem:CN}).\footnote{Alon's seminal work \cite{alon:99} contains many applications of this theorem, and these have served as a huge inspiration to many, including ourselves.  This inspiration continues into this present paper.}
  
 Let $m_3(n)$ denote the answer to Gardner's \nothree queens version problem on an $n\times n$ chessboard.  Gardner knew the precise value of $m_3(n)$ for small values of $n$, and a few more precise values were given in \cite{cooperpikhurkoschmittwarrington:14}.  Subsequently, Rob Pratt \cite{pratt:14} framed the problem as an integer linear programming problem and Don Knuth \cite{knuth:14} brought the power of SAT-solvers to bear on the problem.  Their results are cumulatively stated as an entry in the On-line Encyclopedia of Integer Sequences \cite[{\fontfamily{ppl}\selectfont A219760}]{sloane}.  All values for which $m_3(n)$ is known precisely are given in Table \ref{table:1}.
 
\begin{table}
  \begin{center}
    \begin{tabular}{cccccccccc}\toprule
      $n$       & 1 & 2 & 3 & 4 & 5 & 6 & 7 & 8 & 9 \\
      $m_3(n)$  & 1 & 4 & 4 & 4 & 6 & 6 & 8 &  9 &  10\\\midrule
      $n$ & 10 &       11 & 12 & 13      & 14 & 15 & 16 & 17 & 18\\
      $m_3(n)$ & 10 & 12 & 12 & 14 & 15& 
      16&17& 18 & 18\\\midrule
      $n$ & 19 & 20 & 21 & 22 & 23 & 24 & 25 & 26 & 27\\
      $m_3(n)$ & 20 & 21 & 22 & 23 &24 & 25 & 26 & 26 & 28\\\bottomrule
    \end{tabular}
  \end{center}
  \caption{All known values of $m_3(n)$.}
  \label{table:1}
\end{table}

The data given in Table \ref{table:1} suggested to us that for $n$ odd and $n \geq 3$ that we should have
$m_3(n) \geq n+1$; Gardner conjectured the same \cite{gardner:archives}\footnote{In a typed letter that Martin Gardner wrote to John H. Conway dated 2 June 1975, he wrote, ``It would be nice if the minimum for all odd $n$ were $n+1$, and for all even $n$, $n$ or $n+2$."}.  We prove the conjecture in the case that $n$ is congruent to 1 modulo 4 as our main result.

\begin{theorem}\label{theorem:1mod4}
For $n$ congruent to 1 modulo 4 and $n \geq 5$, we have $m_3(n) \geq n+1$.
\end{theorem} 

The proof of our main result relies on a mixture of the polynomial method, combinatorial arguments and linear algebra, some of which are similar to those found in \cite{cooperpikhurkoschmittwarrington:14}.  However, a new, key ingredient is a recent strengthening of Alon's Combinatorial Nullstellensatz due to Bogdan Nica \cite{nica:22}.  Before stating Nica's theorem, we recall Part II of the Combinatorial Nullstellensatz (i.e., the Non-vanishing Corollary).

\begin{theorem}\label{theorem:CN}[Combinatorial Nullstellensatz, Theorem 1.2~\cite{alon:99}]
  Let $F$ be an arbitrary field, and let $f = f ( x_1 , \ldots , x_n
  )$ be a polynomial in $F [ x_1 , \ldots , x_n ]$.  Suppose the
  degree $\deg (f)$ of $f$ is $\sum_{i = 1}^n t_i$, where each $t_i$
  is a non-negative integer, and suppose the coefficient of $\prod_{i =
    1}^n x_i^{t_i}$ in $f$ is nonzero.  Then, if $S_1 , \ldots , S_n$
  are subsets of $F$ with $|S_i| > t_i$, there are $s_1 \in S_1 ,
  \ldots , s_n \in S_n$ so that $f ( s_1 , \ldots , s_n ) \neq 0$.
\end{theorem}

We say that a subset $S$ of a field $F$ is a {\it zero-sum} if the sum of the elements in $S$ is the zero of the field.  In the case that each of the $S_i$ is a zero-sum subset of the field, Nica \cite{nica:22} showed that the same conclusion in Theorem \ref{theorem:CN} holds under a slightly weaker degree restriction placed upon the degree of the polynomial.  Nica's result is as follows.

\begin{theorem}\label{theorem:CNZS}[Combinatorial Nullstellensatz for Zero-sum Grids, \cite{nica:22}]
  Let $F$ be an arbitrary field, and let $f = f ( x_1 , \ldots , x_n
  )$ be a polynomial in $F [ x_1 , \ldots , x_n ]$.  Suppose the
  degree $\deg (f)$ of $f$ is $1+\sum_{i = 1}^n t_i$, where each $t_i$
  is a non-negative integer, and suppose the coefficient of $\prod_{i =
    1}^n x_i^{t_i}$ in $f$ is nonzero.  Then, if $S_1 , \ldots , S_n$
  are zero-sum subsets of $F$ with $|S_i| > t_i$, there are $s_1 \in S_1 ,
  \ldots , s_n \in S_n$ so that $f ( s_1 , \ldots , s_n ) \neq 0$.
\end{theorem}

This Nullstellensatz for Zero-sum Grids is, in fact, a particular case of a more expansive theorem that we won't discuss here.  Roughly, it says that we may relax the degree constraints on the polynomial to reach the same conclusion of the Combinatorial Nullstellensatz whenever the grid is ``structured", see \cite{nica:22}.  As far as we know, the proof of Theorem \ref{theorem:1mod4} represents the first application of this particular generalization of the Combinatorial Nullstellensatz \cite{nica:22} outside of that paper.  It is interesting to us whenever a generalization of Alon's Nullstellensatz is truly needed for a combinatorial problem, as is the case for our main result.

This paper is organized as follows.  In Section \ref{section:definitions} we give the necessary definitions and a sketch of the proof of Theorem \ref{theorem:1mod4}.  In Section \ref{section:maintheorem} we prove Theorem \ref{theorem:1mod4}.  In Section \ref{section:conclusion} we provide some further remarks and also include insights into the difficulty of showing $m_3(n) \geq n+1$ for the case when $n$ is 3 modulo 4.
\section{Definitions and notation}\label{section:definitions}

\begin{figure}\label{figure:lonely5}
    \centering
    \input{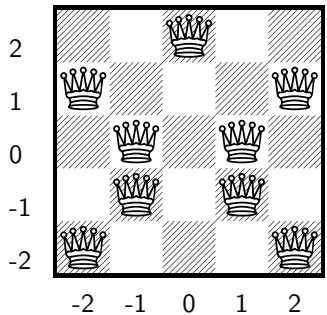}
    \caption{An example of a good placement with one lonely queen on a $5 \times 5$ board.}
\end{figure}

We adopt many of the definitions and notation originally used in \cite{cooperpikhurkoschmittwarrington:14}, with one important exception.  We consider the infinite square ${\mathbb Z}$-lattice as a {\it chessboard} and its vertices as {\it squares} of the chessboard.  A {\it board} $B$ is a finite subset of the chessboard.  For $n$ odd, let $\ZB{n}$ denote the board $[-\frac{n-1}{2},\frac{n-1}{2}] \times [-\frac{n-1}{2},\frac{n-1}{2}].${\footnote {This is the key difference to the set-up as compared to that of \cite{cooperpikhurkoschmittwarrington:14} as we now are considering a zero-sum grid, thus facilitating the use of Theorem \ref{theorem:CNZS}.}}  We may refer to a square of $\ZB{n}$ by the coordinates $(x,y)$ of its corresponding vertex.   As we are interested in the queens version of the problem, the lines that we concern ourselves with have slope $0, +1, -1,$ or $\infty$ and contain vertices of the lattice --- so, throughout we use {\it line} to refer to a line of this type.  Any subset $S$ of the infinite square lattice may be considered a {\it placement of queens}, or {\it placement} for short, by imagining a queen on each corresponding square of the chessboard.  The {\it size} of a placement $S$ is its cardinality $|S|$.  We say that two queens of a placement ${\cal Q}$ {\it define} a line if they lie on the same row, column or diagonal.  In such a way, the placement ${\cal Q}$ {\it defines} a set of lines, the set of lines defined by the collinear pairs of ${\cal Q}$.  A line is said to {\it cover} a square if the coordinates of the square is in the zero locus of the line.  A placement is called  {\it good} if does not contain $3$ queens in a line and loses this property upon the addition of a queen to an unoccupied square.  A queen is said to be {\it lonely} if she is not collinear with any other queen.  We illustrate a placement with a lonely queen in Figure \ref{figure:lonely5}, where there exists one lonely queen in the topmost row (i.e., she has coordinates $(0,2)$).

We now sketch the proof of Theorem \ref{theorem:1mod4}.  We will consider a good placement ${\cal Q}$ of size at most $4k+1$.  The proof will divide into two cases, which roughly depend on the number of lonely queens.  In the first case, when the number of lonely queens is not one, when $q < 4k+1$ or when the number of lines defined by ${\cal Q}$ is not maximized, ${\cal Q}$ will define a set of lines and together these lines (and perhaps some additional lines) will be used to construct a polynomial that vanishes on the squares of the chessboard.  However, we will be able to show that the coefficient on a leading monomial of this polynomial, which is of `small degree,' is non-zero, thus obtaining a contradiction to Alon's Combinatorial Nullstellensatz.  (The proof in the first case is essentially that found in \cite{cooperpikhurkoschmittwarrington:14} for the proof of Theorem \ref{theorem:CPSW14}.) The second case, when the number of lonely queens is one, $q=4k+1$ and the number of lines defined by ${\cal Q}$ is maximized, is more involved.  In this case, Alon's Nullstellensatz {\it will be insufficient}.  Like in the previous case, we use the set of lines defined by ${\cal Q}$ to construct several polynomials. For each of the four possible slopes of a line passing through the lonely queen, we define a polynomial so that it vanishes on all the squares of the chessboard.  For each polynomial, we calculate the coefficient of an appropriate monomial, which is in terms of the coordinates of the lonely queen and sums of the intercepts of the defining lines of ${\cal Q}$.  To compute this coefficient requires nothing more sophisticated than the Binomial Theorem.  If one of these four coefficients is non-zero, then we obtain a contradiction to Nica's Combinatorial Nullstellensatz for Zero-sum Grids.  If each of these four coefficients is zero, then, together with some geometric and combinatorial arguments that yield additional equations, we will be able to construct a homogeneous system of linear equations whose solution yields the location of the lonely queen to be centered on the board (i.e., she has coordinates $(0,0)$) along with some other attributes of the placement.  This structural information allows for a combinatorial argument that finishes the proof.

\section{Proof of Main Theorem}\label{section:maintheorem}

{\sc Proof of Theorem~\ref{theorem:1mod4}}


Let $n=4k+1$ and $k \geq 1$.  Let ${\cal Q}$ be a good placement on $\ZB{n}$ with size $q = |{\cal Q}| \leq 4k+1$.  Let $\QQ'$
denote the (possibly empty) subset of lonely queens in ${\cal Q}$.  Let $|\QQ'|=q'$.  We work towards a contradiction.

{\bf Case 1:} $q' \neq 1$, $q < 4k+1$, or $q'=1$ and some queen(s) that is not lonely participates in defining fewer than $4$ lines.

We proceed by constructing a polynomial $f(x,y)$ of total degree $8k$ that vanishes on each square $(x,y)\in \ZB{n}$.  We will then obtain a contradiction through a suitable application of the Combinatorial Nullstellensatz. 

We construct $f$ as a product of linear factors of three different types.  The {\it first type} consists of the set of lines defined by $\QQ$.  Since the placement $\QQ$ is good, every unoccupied square of $\ZB{n}$ is in the zero locus of at least one line of the first type.

There may be some lonely queens in $\QQ$, which by definition are not on any defining line.    For each $Q_i \in \QQ'$ we will artificially define a new line (of the {\it second type}) that passes through the square occupied by $Q_i$.  While we are free to choose any one of the four possible slopes for each line, it is most convenient to distribute the slopes as evenly as possible. Every occupied square is in the zero locus of at least one line of either of the first two types.

For each of the four possible slopes there are at most $\left\lfloor \frac{q-q'}{2} \right\rfloor$ lines of that slope of the first type and at most $\left\lceil \frac{q'}{4} \right\rceil$ lines of that slope of the second type.  If $q' \neq 1$ or if $q < 4k+1$, these quantities sum to at most $2k$.  If $q'=1$ and some queen(s) that is not lonely participates in defining fewer than $4$ lines, then one of the four possible slope directions has at most $2k-1$ lines and we define a new line through the lonely queen to have that slope.  As necessary, define new, distinct lines (of the {\it third type}) in each of the four directions so that there are exactly $2k$ lines of each slope among the three types.  (The lines of the third type serve only to facilitate the application of the Combinatorial Nullstellensatz; it is immaterial which squares they vanish on.)

Let $\LL = \{L_1,\ldots, L_{8k}\}$ be the set of $8k$ lines and let
$l_i=0$ be the equation in variables $x$ and $y$ defining $L_i$.  We
then define
\begin{equation*}
  f(x,y)  =  \prod_{i=1}^{8k} l_i \in \mathbb{R}[x,y].
\end{equation*}
As desired, the polynomial $f(x,y) = 0$ for every $(x,y)\in \ZB{n}$ as
every unoccupied square is on a line of the first type and every
occupied square is on a line of either the first or second type.  By
construction, the total degree of $f$ is $8k$.  If we group the
factors in $f$ according to slope, we see that $f$ can be rewritten as
\begin{equation*}\label{eq:factor}
  f(x,y) = \prod_{j=1}^{2k} (x-\alpha_j)(y-\beta_j)(x-y-\gamma_j)(x+y-\delta_j)
\end{equation*}
for suitable constants $\alpha_j,\beta_j,\gamma_j,\delta_j$.  The leading monomials of $f(x,y)$ are the same as that of $g(x,y)=x^{2k}y^{2k}(x^2-y^2)^{2k}$.  From the Binomial Theorem, we conclude that the coefficient of the top-degree term $x^{4k}y^{4k}$ is $\pm {2k \choose k}$, i.e. the coefficient is nonzero.

We now apply Theorem~\ref{theorem:CN} to $f(x,y)$, where $t_1=t_2=4k$
and $S_1=S_2=\{-2k,\ldots, 2k\}$, to obtain that there are $s_1 \in S_1, s_2 \in
S_2$ such that $f(s_1,s_2) \neq 0$.  We have reached a contradiction.

{\bf Case 2:} $q=4k+1$, $q'=1$ and all other queens participate in defining a line of each of the four possible slopes.

As we did in Case 1, we employ the polynomial method; this time we will apply Theorem \ref{theorem:CNZS}.

We start by constructing a polynomial $f_1:=f_1(x,y)$ of total degree $8k+1$ that vanishes on each square $(x,y)\in \ZB{n}$.  The polynomial $f_1$ will be a product of linear factors of two different types.  The first type is given by the set of lines defined by $\QQ$.  Since the placement $\QQ$ is good, every unoccupied square of $\ZB{n}$ is in the zero locus of at least one line of the first type.  As given by the restrictions of the case, there is one lonely queen in $\QQ'$ which is not on any defining line.    For this one lonely queen $Q \in \QQ'$ we artificially define a new line that passes through the square occupied by $Q$.  While we are free to choose any one of the four possible slopes for this one line, at this point in the proof we choose slope $\infty$.  Every occupied square is covered by at least one line of either of the two types.

For each of the four possible slopes there are exactly $\frac{4k}{2}=2k$ lines of that slope of the first type and precisely one line of slope $\infty$ of the second type.  (Unlike in Case 1, we need not define lines of the third type.)  

Let $\LL = \{L_1,\ldots, L_{8k+1}\}$ be our set of $8k+1$ lines and let
$l_i=0$ be the equation in variables $x$ and $y$ defining $L_i$.  We
then define
\begin{equation*}
  f_1(x,y)  =  \prod_{i=1}^{8k+1} l_i \in \mathbb{R}[x,y].
\end{equation*}
As desired, the polynomial $f_1(x,y) = 0$ for every $(x,y)\in \ZB{n}$ as
every unoccupied square is on a line of the first type and every
occupied square is on a line of either the first or second type.  By
construction, the total degree of $f$ is $8k+1$.  Let $(\alpha_0,\beta_0)$ denote the square occupied by the lonely queen.  If we group the
factors in $f_1$ according to slope, we see that $f_1$ can be rewritten as
\begin{equation*}\label{eq:f1factor}
  f_1(x,y) = (x-\alpha_0)\prod_{j=1}^{2k} (x-\alpha_j)(y-\beta_j)(x-y-\gamma_j)(x+y-\delta_j)
\end{equation*}
for suitable constants $\alpha_j,\beta_j,\gamma_j,\delta_j$.

If can we conclude that
the coefficient of the term $x^{4k}y^{4k}$ is nonzero, then we may apply Theorem~\ref{theorem:CNZS} to $f_1(x,y)$, where $t_1=t_2=4k$
and $S_1=S_2=\{-2k,\ldots, 2k\}$ (zero-sum sets), to obtain that there are $s_1 \in S_1, s_2 \in
S_2$ such that $f_1(s_1,s_2) \neq 0$, thus obtaining a contradiction.

We thus `expand' $f_1(x,y)$ and `collect' like terms so that we might determine the coefficient of $x^{4k}y^{4k}$.  Since we are only interested in the coefficient on this term, we focus our analysis only on it.  Note that this term has degree one less than the degree of $f_1$.  So, to obtain such a monomial in the expansion, from the $8k+1$ linear factors, we must choose the $x$-variable $4k$ times, the $y$-variable $4k$ times and thus some constant once.  We think of choosing that constant first and so partition our analysis based upon whether the constant that we have chosen is some $\alpha_j, \beta_j, \gamma_j$ or $\delta_j$.

\begin{enumerate}
\item Choose some $\alpha_j$ first for some $0 \leq j \leq 2k$.\\
The remaining factors with their constant terms removed (since we can't choose them) are:  $$x^{2k}y^{2k}(x-y)^{2k}(x+y)^{2k}.$$
This equals $$x^{2k}y^{2k}(x^2-y^2)^{2k}=x^{2k}y^{2k}\sum_{\ell=0}^{2k}{2k \choose \ell}(x^2)^{\ell}(-y^2)^{2k-\ell}.$$
The only choice of $\ell$ which yields the desired monomial is $\ell =k$, which gives a coefficient of ${2k \choose k}(-1)^k$.  This is the contribution for each $\alpha_j$.  Thus, we have a total contribution to the coefficient of the desired monomial of $(-1)^k{2k \choose k}\sum_{j=0}^{2k}-\alpha_j.$

\item Choose some $\beta_j$ first for some $1 \leq j \leq 2k$.\\
The remaining factors with their constant terms removed (since we can't choose them) are:  $$x^{2k+1}y^{2k-1}(x-y)^{2k}(x+y)^{2k}.$$
This equals $$x^{2k+1}y^{2k-1}(x^2-y^2)^{2k}=x^{2k+1}y^{2k-1}\sum_{\ell=0}^{2k}{2k \choose \ell}(x^2)^{\ell}(-y^2)^{2k-\ell}.$$
There is no choice of $\ell$ which yields the desired monomial, which is true for each $\beta_j$.  Thus, we have a total contribution to the coefficient of the desired monomial of $0.$

\item Choose some $\gamma_j$ first for some $1 \leq j \leq 2k$.\\

The remaining factors with their constant terms removed (since we can't choose them) are:  $$x^{2k+1}y^{2k}(x-y)^{2k-1}(x+y)^{2k}.$$
This equals
\begin{eqnarray*}
x^{2k+1}y^{2k}(x+y)(x^2-y^2)^{2k-1} &= &x^{2k+1}y^{2k}(x+y)\sum_{\ell=0}^{2k-1}{2k-1 \choose \ell}(x^2)^{\ell}(-y^2)^{2k-1-\ell}\\
   & = &x^{2k+2}y^{2k}\sum_{\ell=0}^{2k-1}{2k-1 \choose \ell}(x^2)^{\ell}(-y^2)^{2k-1-\ell} \\ & & +x^{2k+1}y^{2k+1}\sum_{\ell=0}^{2k-1}{2k-1 \choose \ell}(x^2)^{\ell}(-y^2)^{2k-1-\ell}.
\end{eqnarray*}
Only in the first of the two summands is there a choice of $\ell$ which yields the desired monomial; it is $\ell =k-1$, which gives a coefficient of ${2k-1 \choose k-1}(-1)^k$.  This is the contribution for each $\gamma_j$.  Thus, we have a total contribution to the coefficient of the desired monomial of $(-1)^k{2k-1 \choose k-1}\sum_{j=1}^{2k}-\gamma_j.$

\item Choose some $\delta_j$ first for some $1 \leq j \leq 2k$.\\

The remaining factors with their constant terms removed (since we can't choose them) are:  $$x^{2k+1}y^{2k}(x-y)^{2k}(x+y)^{2k-1}.$$
This equals
\begin{eqnarray*}
x^{2k+1}y^{2k}(x-y)(x^2-y^2)^{2k-1} &= &x^{2k+1}y^{2k}(x-y)\sum_{\ell=0}^{2k-1}{2k-1 \choose \ell}(x^2)^{\ell}(-y^2)^{2k-1-\ell} \\ 
   & = & x^{2k+2}y^{2k}\sum_{\ell=0}^{2k-1}{2k-1 \choose \ell}(x^2)^{\ell}(-y^2)^{2k-1-\ell}\\
   & & - x^{2k+1}y^{2k+1}\sum_{\ell=0}^{2k-1}{2k-1 \choose \ell}(x^2)^{\ell}(-y^2)^{2k-1-\ell}.
\end{eqnarray*}

Only in the first of the two summands is there a choice of $\ell$ which yields the desired monomial; it is $\ell =k-1$, which gives a coefficient of ${2k-1 \choose k-1}(-1)^k$.  This is the contribution for each $\delta_j$.  Thus, we have a total contribution to the coefficient of the desired monomial of $(-1)^k{2k-1 \choose k-1}\sum_{j=1}^{2k}-\delta_j.$

\end{enumerate}

The sum of these four contributions is the coefficient $\coef{f_1}{(4k,4k)}$ on $x^{4k}y^{4k}$ in $f_1(x,y)$.  It is 
\begin{eqnarray*}\label{f1coefficient}
\coef{f_1}{(4k,4k)} &=& 
(-1)^k{2k \choose k}\sum_{j=0}^{2k}-\alpha_j+(-1)^k{2k-1 \choose k-1}\sum_{j=1}^{2k}-\gamma_j+(-1)^k{2k-1 \choose k-1}\sum_{j=1}^{2k}-\delta_j.
\end{eqnarray*}

If this coefficient is non-zero, then by Theorem \ref{theorem:CNZS} we are done.  Thus, assume it is zero.

We repeat the above procedure for the placement $\QQ$ with the one lonely queen at $(\alpha_0, \beta_0)$ by defining a similar polynomial $f_2(x,y)$, this time with a line of slope $0$ through the lonely queen, as follows.

Let \begin{equation*}\label{eq:f2factor}
  f_2(x,y) = (y-\beta_0)\prod_{j=1}^{2k} (x-\alpha_j)(y-\beta_j)(x-y-\gamma_j)(x+y-\delta_j)
\end{equation*}
for suitable constants $\alpha_j,\beta_j,\gamma_j,\delta_j$.

By repeating the above procedure or, perhaps more simply, by noting the symmetries between $f_1(x,y)$ and $f_2(x,y)$, the coefficient $\coef{f_2}{(4k,4k)}$ on $x^{4k}y^{4k}$ in $f_2(x,y)$ is 
\begin{eqnarray*}\label{f2coefficient}
\coef{f_2}{(4k,4k)} & = & (-1)^k{2k \choose k}\sum_{j=0}^{2k}-\beta_j+(-1)^k{2k-1 \choose k-1}\sum_{j=1}^{2k}-\gamma_j+(-1)^k{2k-1 \choose k-1}\sum_{j=1}^{2k}-\delta_j.
\end{eqnarray*}

If this coefficient is non-zero, then by Theorem \ref{theorem:CNZS} we are done.  Thus, assume it is zero.



We repeat the above procedure for the placement $\QQ$ with the one lonely queen at $(\alpha_0, \beta_0)$ by defining a similar polynomial $f_3(x,y)$, this time with a line of slope $+1$ through the lonely queen, as follows.

Let \begin{equation*}\label{eq:f3factor}
  f_3(x,y) = (x-y-\gamma_0)\prod_{j=1}^{2k} (x-\alpha_j)(y-\beta_j)(x-y-\gamma_j)(x+y-\delta_j)
\end{equation*}
for suitable constants $\alpha_j,\beta_j,\gamma_j,\delta_j$.

As above, we `expand' $f_3(x,y)$ and `collect' like terms so that we might determine the coefficient of $x^{4k}y^{4k}$.  Since we are only interested in the coefficient on this term, we focus our analysis only on it.  Note that this term has degree one less than the degree of $f_3$.  So, to obtain such a monomial in the expansion, from the $8k+1$ linear factors, we must choose the $x$-variable $4k$ times, the $y$-variable $4k$ times and thus some constant once.  We think of choosing that constant first and so partition our analysis based upon whether the constant that we have chosen is some $\alpha_j, \beta_j, \gamma_j$ or $\delta_j$.

\begin{enumerate}
\item Choose some $\alpha_j$ first for some $1 \leq j \leq 2k$.\\
The remaining factors with their constant terms removed (since we can't choose them) are:  $$x^{2k-1}y^{2k}(x-y)^{2k+1}(x+y)^{2k}.$$
This equals 
\begin{eqnarray*}
x^{2k-1}y^{2k}(x-y)(x^2-y^2)^{2k} & = &x^{2k-1}y^{2k}(x-y)\sum_{\ell=0}^{2k}{2k \choose \ell}(x^2)^{\ell}(-y^2)^{2k-\ell}\\
 & =& x^{2k}y^{2k}\sum_{\ell=0}^{2k}{2k \choose \ell}(x^2)^{\ell}(-y^2)^{2k-\ell} \\
 & & - x^{2k-1}y^{2k+1}\sum_{\ell=0}^{2k}{2k \choose \ell}(x^2)^{\ell}(-y^2)^{2k-\ell}.
\end{eqnarray*}

Only in the first of the two summands is there a choice of $\ell$ which yields the desired monomial; it is $\ell =k$, which gives a coefficient of ${2k \choose k}(-1)^k$.  This is the contribution for each $\alpha_j$.  Thus, we have a total contribution to the coefficient of the desired monomial of $(-1)^k{2k \choose k}\sum_{j=1}^{2k}-\alpha_j.$

\item Choose some $\beta_j$ first for some $1 \leq j \leq 2k$.\\
The remaining factors with their constant terms removed (since we can't choose them) are:  $$x^{2k}y^{2k-1}(x-y)^{2k+1}(x+y)^{2k}.$$
This equals 
\begin{eqnarray*}
x^{2k}y^{2k-1}(x-y)(x^2-y^2)^{2k} & = &x^{2k}y^{2k-1}(x-y)\sum_{\ell=0}^{2k}{2k \choose \ell}(x^2)^{\ell}(-y^2)^{2k-\ell}\\
 & =& x^{2k+1}y^{2k-1}\sum_{\ell=0}^{2k}{2k \choose \ell}(x^2)^{\ell}(-y^2)^{2k-\ell} \\
 & & - x^{2k}y^{2k}\sum_{\ell=0}^{2k}{2k \choose \ell}(x^2)^{\ell}(-y^2)^{2k-\ell}.
\end{eqnarray*}

Only in the second of the two summands is there a choice of $\ell$ which yields the desired monomial; it is $\ell =k$, which gives a coefficient of $(-1){2k \choose k}(-1)^k$.  This is the contribution for each $\beta_j$.  Thus, we have a total contribution to the coefficient of the desired monomial of $(-1)^{k+1}{2k \choose k}\sum_{j=1}^{2k}-\beta_j.$

\item Choose some $\gamma_j$ first for some $0 \leq j \leq 2k$.\\

The remaining factors with their constant terms removed (since we can't choose them) are:  $$x^{2k}y^{2k}(x-y)^{2k}(x+y)^{2k}.$$
This equals
\begin{eqnarray*}
x^{2k}y^{2k}(x^2-y^2)^{2k} &= &x^{2k}y^{2k}\sum_{\ell=0}^{2k}{2k \choose \ell}(x^2)^{\ell}(-y^2)^{2k-\ell}.
\end{eqnarray*}

The only choice of $\ell$ which yields the desired monomial is $\ell =k$, which gives a coefficient of ${2k \choose k}(-1)^k$.  This is the contribution for each $\gamma_j$.  Thus, we have a total contribution to the coefficient of the desired monomial of $(-1)^k{2k \choose k}\sum_{j=0}^{2k}-\gamma_j.$

\item Choose some $\delta_j$ first for some $1 \leq j \leq 2k$.\\

The remaining factors with their constant terms removed (since we can't choose them) are:  $$x^{2k}y^{2k}(x-y)^{2k+1}(x+y)^{2k-1}.$$

This equals
\begin{eqnarray*}
x^{2k}y^{2k}(x-y)^2(x^2-y^2)^{2k-1} &= &x^{2k}y^{2k}(x^2-2xy+y^2)\sum_{\ell=0}^{2k-1}{2k-1 \choose \ell}(x^2)^{\ell}(-y^2)^{2k-1-\ell}\\
& = & x^{2k+2}y^{2k}\sum_{\ell=0}^{2k-1}{2k-1 \choose \ell}(x^2)^{\ell}(-y^2)^{2k-1-\ell} \\
& & -2x^{2k+1}y^{2k+1}\sum_{\ell=0}^{2k-1}{2k-1 \choose \ell}(x^2)^{\ell}(-y^2)^{2k-1-\ell}\\
& & +x^{2k}y^{2k+2}\sum_{\ell=0}^{2k-1}{2k-1 \choose \ell}(x^2)^{\ell}(-y^2)^{2k-1-\ell}.
\end{eqnarray*}

Amongst the three summands, there is a choice of $\ell=k-1$ in the first, no choice in the second and a choice of $\ell=k$ in the third.  This yields a coefficient of $(-1)^k{2k-1 \choose k-1}+(-1)^{k-1}{2k-1 \choose k}=0.$  This is the contribution for each $\delta_j$.  Thus, we have a total contribution to the coefficient of the desired monomial of $0.$

\end{enumerate}

The sum of these four contributions is the coefficient $\coef{f_3}{(4k,4k)}$ on $x^{4k}y^{4k}$ in $f_3(x,y)$.  It is 
\begin{eqnarray*}\label{f3coefficient}
\coef{f_3}{(4k,4k)} &=& 
(-1)^k{2k \choose k}\sum_{j=1}^{2k}-\alpha_j+(-1)^{k+1}{2k \choose k}\sum_{j=1}^{2k}-\beta_j+(-1)^k{2k \choose k}\sum_{j=0}^{2k}-\gamma_j.
\end{eqnarray*}

If this coefficient is non-zero, then by Theorem \ref{theorem:CNZS} we are done.  Thus, assume it is zero.

We repeat the above procedure for the placement $\QQ$ with the one lonely queen at $(\alpha_0, \beta_0)$ by defining a similar polynomial $f_4(x,y)$, this time with a line of slope $-1$ through the lonely queen, as follows.

Let \begin{equation*}\label{eq:f4factor}
  f_4(x,y) = (x+y-\delta_0)\prod_{j=1}^{2k} (x-\alpha_j)(y-\beta_j)(x-y-\gamma_j)(x+y-\delta_j)
\end{equation*}
for suitable constants $\alpha_j,\beta_j,\gamma_j,\delta_j$.


By repeating the above procedure or, perhaps more simply, by noting the symmetries between $f_3(x,y)$ and $f_4(x,y)$, the coefficient $\coef{f_4}{(4k,4k)}$ on $x^{4k}y^{4k}$ in $f_4(x,y)$ is 
\begin{eqnarray*}\label{f4coefficient}
\coef{f_4}{(4k,4k)} &=& 
(-1)^k{2k \choose k}\sum_{j=1}^{2k}-\alpha_j+(-1)^{k}{2k \choose k}\sum_{j=1}^{2k}-\beta_j+(-1)^k{2k \choose k}\sum_{j=0}^{2k}-\delta_j.
\end{eqnarray*}

If this coefficient is non-zero, then by Theorem \ref{theorem:CNZS} we are done.  Thus, assume it is zero.

We have now reached the following system of linear equations: 
\begin{eqnarray}\label{collectedcoefficients}
\coef{f_1}{(4k,4k)}=0, \coef{f_2}{(4k,4k)}=0, \coef{f_3}{(4k,4k)}=0, \coef{f_4}{(4k,4k)}=0.
\end{eqnarray}

Next we generate four additional linear equations via some geometric and combinatorial observations.

Consider the lonely queen located on square $(\alpha_0,\beta_0)$: the values of $\alpha_0$ and $\beta_0$ determine the values of $\gamma_0$ and $\delta_0$ as follows.  The line of slope $+1$ that goes through the square $(\alpha_0, \beta_0)$ has equation $$y-\beta_0=1(x-\alpha_0), ~~~x-y-(\alpha_0-\beta_0)=0$$ and the line of slope $-1$ that goes through the square $(\alpha_0,\beta_0)$ has equation $$y-\beta_0=-1(x-\alpha_0),~~~x+y-(\alpha_0+\beta_0)=0.$$
As a result, we have 

\begin{eqnarray}
\alpha_0 - \beta_0-\gamma_0=0, \label{lonelylocation1} \\
\alpha_0 + \beta_0-\delta_0=0. \label{lonelylocation2}
\end{eqnarray}

Now consider the other $4k$ queens of ${\cal Q} \setminus {\cal Q'}$ (i.e., those that are not lonely).  For each such queen there exists an $\alpha \in \{\alpha_1, \ldots , \alpha_{2k}\}$ and $\beta \in \{\beta_1, \ldots , \beta_{2k}\}$ that give her coordinates.  The line of slope $+1$ that goes through the square $(\alpha, \beta)$ has equation $$y-\beta=1(x-\alpha), ~~~x-y-(\alpha-\beta)=0$$ and the line of slope $-1$ that goes through the square $(\alpha,\beta)$ has equation $$y-\beta=-1(x-\alpha),~~~x+y-(\alpha+\beta)=0.$$
As a result, we have $\gamma=\alpha - \beta$ for some $\gamma \in \{\gamma_1, \ldots ,\gamma_{2k}\}$ and $\delta= \alpha + \beta$ for some $\delta \in \{\delta_1, \ldots, \delta_{2k}\}$.  As each such diagonal line is defined by two queens, upon considering the $4k$ equations deriving from the $-1$-slope lines each element of $\{\gamma_1, \ldots ,\gamma_{2k}\}$ occurs twice in this set of equations; similarly, in the $4k$ equations deriving from the $+1$-slope lines each element of $\{\delta_1, \ldots, \delta_{2k}\}$ occurs twice.  Thus, we can write the following:

\begin{eqnarray}\label{eqn:slopesminus1}
 \sum_{Q \in {\cal Q}\setminus {\cal Q'}}(\alpha-\beta) = 2 \sum_{i=1}^{2k}\gamma_i,
 \end{eqnarray}
and
\begin{eqnarray}\label{eqn:slopesplus1}
   \sum_{Q \in {\cal Q}\setminus {\cal Q'}}(\alpha+\beta) = 2 \sum_{i=1}^{2k}\delta_i. 
\end{eqnarray}

The restrictions of Case 2 give that each queen in ${\cal Q}\setminus{\cal Q'}$ is contained in both a vertical line and a horizontal line.  As a result, each $\alpha \in \{\alpha_1, \ldots , \alpha_{2k}\}$ and each $\beta \in \{\beta_1, \ldots , \beta_{2k}\}$ appears twice on the left side of each of Equation \ref{eqn:slopesminus1} and Equation \ref{eqn:slopesplus1}.  This enables us to rewrite the left side of Equations \ref{eqn:slopesminus1} and \ref{eqn:slopesplus1} to obtain

\begin{eqnarray}\label{eqn:slopesminus1A}
2\sum_{i=1}^{2k}\alpha_i - 2\sum_{i=1}^{2k} \beta_i = 2 \sum_{i=1}^{2k}\gamma_i,
\end{eqnarray}
and
\begin{eqnarray}\label{eqn:slopesplus1A}
  2\sum_{i=1}^{2k}\alpha_i + 2\sum_{i=1}^{2k} \beta_i = 2 \sum_{i=1}^{2k}\delta_i.  
\end{eqnarray}

We may scale and rewrite Equations \ref{eqn:slopesminus1A}-\ref{eqn:slopesplus1A} as

\begin{eqnarray}\label{eqn:slopesminus1B}
 \sum_{i=1}^{2k}\alpha_i - \sum_{i=1}^{2k} \beta_i -  \sum_{i=1}^{2k}\gamma_i =0,   
\end{eqnarray}

and
\begin{eqnarray}\label{eqn:slopesplus1B}
   \sum_{i=1}^{2k}\alpha_i + \sum_{i=1}^{2k} \beta_i -  \sum_{i=1}^{2k}\delta_i =0. 
\end{eqnarray}


At this point, we have now generated a homogeneous system of eight linear equations - these are Equations \ref{collectedcoefficients}, \ref{lonelylocation1}, \ref{lonelylocation2}, \ref{eqn:slopesminus1B}, \ref{eqn:slopesplus1B} - in the variables $\alpha_0, \beta_0, \gamma_0, \delta _0, \sum_{j=1}^{2k}\alpha_j, \sum_{j=1}^{2k}\beta_j, \sum_{j=1}^{2k}\gamma_j, \sum_{j=1}^{2k}\delta_j.$  For notational convenience we set $\omega := (-1)^k{2k \choose k}$ and express these equations using the following augmented matrix.

\[
  \mathbf{[{\bf A}~~0]} = 
    \begin{blockarray}{ccccccccc}
         \alpha_0 & \beta_0 & \gamma_0 & \delta _0 & \sum_{j=1}^{2k}\alpha_j &  \sum_{j=1}^{2k}\beta_j  &  \sum_{j=1}^{2k}\gamma_j  &  \sum_{j=1}^{2k}\delta_j  &   \\
      \begin{block}{[ccccccccc]}
    \omega & 0 & 0 & 0 & \omega & 0 & \omega/2 & \omega/2 & 0 \\
    0 & \omega & 0 & 0 & 0 & \omega & \omega/2 & \omega/2 & 0 \\
    0 & 0 & \omega & 0 & \omega & -\omega & \omega & 0 & 0 \\
    0 & 0 & 0 & \omega & \omega & \omega & 0 & \omega & 0 \\
    1 & -1 & -1 & 0 & 0 & 0 & 0 & 0 & 0 \\
    1 & 1 & 0 & -1 & 0 & 0 & 0 & 0 & 0\\
    0 & 0 & 0 & 0 & 1 & -1 & -1 & 0 & 0 \\
    0 & 0 & 0 & 0 & 1 & 1 & 0 & -1 & 0 \\
      \end{block}
    \end{blockarray}
\]

Guassian elimination yields the following row reduced $8 \times 8$ coefficient matrix:

\[
  \mathbf{{\bf A}} \sim 
    \begin{blockarray}{cccccccc}
         \alpha_0 & \beta_0 & \gamma_0 & \delta _0 & \sum_{j=1}^{2k}\alpha_j &  \sum_{j=1}^{2k}\beta_j  &  \sum_{j=1}^{2k}\gamma_j  &  \sum_{j=1}^{2k}\delta_j  \\
      \begin{block}{[cccccccc]}
      1 & 0 & 0 & 0 & 0 & 0 & 0 & 1 \\
      0 & 1 & 0 & 0 & 0 & 0 & 0 & 1 \\
      0 & 0 & 1 & 0 & 0 & 0 & 0 & 0 \\
      0 & 0 & 0 & 1 & 0 & 0 & 0 & 2 \\
      0 & 0 & 0 & 0 & 1 & 0 & 0 & -0.5\\
      0 & 0 & 0 & 0 & 0 & 1 & 0 & -0.5\\
      0 & 0 & 0 & 0 & 0 & 0 & 1 & 0\\
      0 & 0 & 0 & 0 & 0 & 0 & 0 & 0\\
      \end{block}
    \end{blockarray}
\]
We see that the null space of $A$ is spanned by the following vector,
\[ \begin{bmatrix}
        1 & 1 & 0 & 2 & -0.5 & -0.5 & 0 & -1
    \end{bmatrix}^{\intercal}.\]

Note, for a vector ${\bf x}$ {\it not} in the null space, at least one entry of $A{\bf x}$ is nonzero. If this nonzero entry is among the first $4$ entries, then at least one of $\coef{f_1}{(4k,4k)}, \coef{f_2}{(4k,4k)}, \coef{f_3}{(4k,4k)}, \coef{f_4}{(4k,4k)}$ is nonzero, a contradiction to Theorem \ref{theorem:CNZS}. If this nonzero entry is among the latter $4$ entries, then it contradicts a relationship we found in one of Equations \ref{lonelylocation1}, \ref{lonelylocation2}, \ref{eqn:slopesminus1B}, \ref{eqn:slopesplus1B}.

The vectors in the null space of $A$ imply that the lonely queen has coordinates $(s, s)$ for $s \in \{-2k,\ldots , 2k\}$, so the lonely queen is on the line $y = x$, i.e. the `back-diagonal' of the chessboard.  If we rotate the placement by 90-degrees counter-clockwise, then the placement we obtain satisfies the conditions of Case 2 and has the lonely queen on the `forward-diagonal', i.e., on the line $y=-x$.  Applying all of the previous arguments to this new placement shows that the lonely queen must be on the back-diagonal.  Thus, the lonely queen is at the center, i.e., she has coordinates $(0,0)$.  This implies that the only vector ${\bf x}$ in the null space of $A$ which we need to be concerned with is the zero vector.  Thus, 

\begin{eqnarray}\label{eqn:nullA}
\alpha_0, \beta_0, \gamma_0, \delta _0, \sum_{j=1}^{2k}\alpha_j, \sum_{j=1}^{2k}\beta_j, \sum_{j=1}^{2k}\gamma_j, \sum_{j=1}^{2k}\delta_j=0.
\end{eqnarray}

\begin{figure}
    \centering
    \input{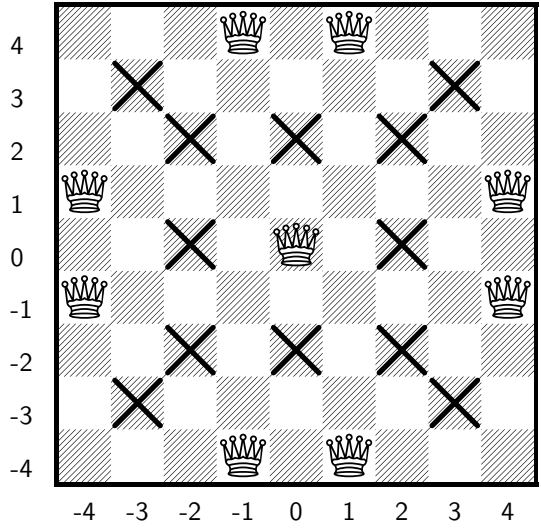}
    \caption{A placement of 10 queens on a $9 \times 9$ board which corresponds to the zero vector in Nul $A$.  Sqaures marked with a cross indicate those squares for which the additional placement of a queen would {\it not} yield three-in-a-line.}
    \label{fig:10on9by9withlonelycentered}
\end{figure}

With the location information provided by Equation \ref{eqn:nullA} and the restrictions of Case 2, we now turn to a combinatorial and geometric argument to finish the proof.\footnote{Figure \ref{fig:10on9by9withlonelycentered} is an example of a placement satisfying Equation \ref{eqn:nullA} yet is not good.}

We distinguish between the lines of slope $0$ or $\infty$ defined by ${\cal Q}$ and those of slope $\pm 1$.  To this end, set $U\subseteq \ZB{n}$ to be the set of squares left uncovered by a line of slope $0$ or $\infty$, where the lonely queen {\it does not} define a line. Notice that $U$ defines a rectangular sub-board.  For any index $i \in \{\frac{-(n-1)}{2}, \ldots, \frac{n-1}{2}\}$ (respectively $j \in \{\frac{-(n-1)}{2}, \ldots, \frac{n-1}{2}\}$) let $C_i = \{(i,\ell)\in U :\, \frac{-(n-1)}{2} \leq \ell \leq \frac{n-1}{2}\}$ (respectively $R_j = \{(\ell,j)\in U :\, \frac{-(n-1)}{2} \leq \ell \leq \frac{n-1}{2}\}$).  The sets $C_i$ and $R_j$ keep track of the squares in $U$ for each column and row, respectively.  Let $a < b$ be the minimum and maximum indices, respectively, for which $C_i \neq \emptyset$.  Define $a' < b'$ analogously for the sets $R_j$.   The number of the $C_i$ and $R_j$ that are nonempty is $n-2k=4k+1-2k=2k+1$.  The $8k$ squares that form the set $C_a \cup C_b \cup R_{a'} \cup R_{b'}$ will be referred to as the {\it perimeter of $U$}.
 
Without loss of generality, we may assume $b - a \geq b' - a'$ as otherwise we may rotate the placement by $90^\circ$; also, note that $b-a \geq b'-a' \geq 2k$.  Consider the case that $b-a > b'-a'$ (i.e. the inequality is strict, $U$ is rectangular but not square).  In this case, there are at least $4k+1$ (and at most $4k+2$) empty squares of $C_a \cup C_b$, each of which must be covered by a line of slope $\pm 1$.  However, each such line can cover at most one of these squares.  As there are $4k$ such lines, we fall short of being able to cover each such square.  So, we may conclude that $b-a=b'-a'$. 

If the lonely queen were to occupy a square in $C_a$ (that is, $a=0$), then all vertical lines defined by ${\cal Q}$ would be to the left of center contradicting that $\sum_{i=1}^{2k}\alpha_i=0$.  A similar contradiction would be reached if the lonely queen were to occupy a square in $C_b$ (that is, $b=0$): all vertical lines defined by ${\cal Q}$ would be to the right of center contradicting that $\sum_{i=1}^{2k}\alpha_i=0$.  Thus, the squares of $C_a \cup C_b$ are empty of queens.  Likewise, the squares of  $R_{a'} \cup R_{b'}$ are empty of queens.  Thus, the perimeter of $U$ is empty of queens.

\begin{figure}
   \centering
    \tikzset{every picture/.style={line width=0.75pt}} 
\resizebox{200pt}{!}{%
\begin{tikzpicture}[x=0.75pt,y=0.75pt,yscale=-1,xscale=1]

\draw   (116,24) -- (527,24) -- (527,435) -- (116,435) -- cycle ;
\draw   (143,77) -- (364,77) -- (364,298) -- (143,298) -- cycle ;
\draw    (135,244) -- (166.08,213.09) -- (316,64) ;
\draw    (143,298) -- (364,77) ;
\draw    (364,298) -- (143,77) ;
\draw    (126,225) -- (307,45) ;
\draw    (194,316) -- (375,136) ;
\draw    (213,326) -- (268.65,270.65) -- (394,146) ;

\draw (264,113.4) node [anchor=north west][inner sep=0.75pt]    {$U$};
\draw (120,185.4) node [anchor=north west][inner sep=0.75pt]    {$C_{a}$};
\draw (244,303.4) node [anchor=north west][inner sep=0.75pt]    {$R_{a'}$};
\draw (372,194.4) node [anchor=north west][inner sep=0.75pt]    {$C_{b}$};
\draw (238,52.4) node [anchor=north west][inner sep=0.75pt]    {$R_{b'} \ $};
\draw (140,440.4) node [anchor=north west][inner sep=0.75pt]    {$a$};
\draw (358,443.4) node [anchor=north west][inner sep=0.75pt]    {$b$};
\draw (89,295.4) node [anchor=north west][inner sep=0.75pt]    {$a'$};
\draw (91,66.4) node [anchor=north west][inner sep=0.75pt]    {$b'$};
\draw (121,306.4) node [anchor=north west][inner sep=0.75pt]    {$( a,\ a')$};
\draw (345,306.4) node [anchor=north west][inner sep=0.75pt]    {$( b,\ a')$};
\draw (349,55.4) node [anchor=north west][inner sep=0.75pt]    {$( b,\ b')$};
\draw (125,55.4) node [anchor=north west][inner sep=0.75pt]    {$( a,\ b')$};
\draw (304,442.4) node [anchor=north west][inner sep=0.75pt]    {$0$};
\draw (90,222.4) node [anchor=north west][inner sep=0.75pt]    {$0$};
\draw (512,442.4) node [anchor=north west][inner sep=0.75pt]    {$2k$};
\draw (103,441.4) node [anchor=north west][inner sep=0.75pt]    {$-2k$};
\draw (91,22.4) node [anchor=north west][inner sep=0.75pt]    {$2k$};
\draw (82,420.4) node [anchor=north west][inner sep=0.75pt]    {$-2k$};

\end{tikzpicture}
}
    \caption{The $8k$ squares that form the perimeter of $U$ and some covering lines.}
    \label{fig:updated_tikz_image}
\end{figure}
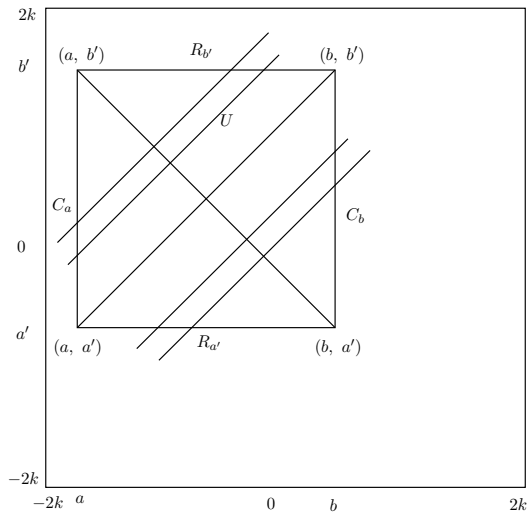

Consider the $8k$ squares that form the perimeter of $U$ as shown in Figure \ref{fig:updated_tikz_image}.  As the perimeter of $U$ is empty of queens, each of these squares must be covered by a line of slope $\pm 1$.  Such a line can cover at most 2 such squares.  With a total of $4k$ diagonal lines ($2k$ of slope $+1$ and $2k$ of slope $-1$) and $8k$ squares, each diagonal line must cover 2 squares.  This forces the $-1$-slope diagonal covering the squares $(a,b')$ and $(b,a')$ and forces the $+1$-slope diagonal covering the squares $(a,a')$ and $(b,b')$ to exist.  Also, a $-1$-slope line that covers a square of $C_a \setminus \{(a,a'),(a,b')\}$ must also cover a square of $R_{a'} \setminus \{(a,a'),(b,a')\}$; likewise,  a $-1$-slope line that covers a square of $C_b \setminus \{(b,a'),(b,b')\}$ must also cover a square of $R_{b'} \setminus \{(a,b'),(b,b')\}$.  Also,  a $+1$-slope line that covers a square of $C_a \setminus \{(a,a'),(a,b')\}$ must also cover a square of $R_{b'} \setminus \{(a,b'),(b,b')\}$; likewise,  a $+1$-slope line that covers a square of $C_b \setminus \{(b,a'),(b,b')\}$ must also cover a square of $R_{a'} \setminus \{(a,a'),(b,a')\}$.  Suppose that there are $p$ lines of slope $-1$ covering squares of $(C_a \setminus \{(a,a'),(a,b')\}) \cup (R_{a'} \setminus \{(a,a'),(b,a')\})$.  The remaining $2k-1-p$ lines of slope $-1$ cover the squares of $(C_b \setminus \{(b,a'),(b,b')\}) \cup (R_{b'} \setminus \{(a,b'),(b,b')\})$, leaving $4k-2-2(2k-1-p)=2p$ squares of this set uncovered, with $p$ squares in $C_b \setminus \{(b,a'),(b,b')\}$ and $p$ squares in $R_{b'} \setminus \{(a,b'),(b,b')\}$.  These squares must be covered by $+1$-slope lines, for an additional $2p$ lines of slope $+1$.  These $2p$ lines together with the $+1$-slope diagonal covering $(a,a')$ and $(b,b')$ give a total of $2p+1$ such lines, which is an odd number of such lines.  However, we said that the number of lines of slope $+1$ is $2k$, which is even, and this contradiction finishes the proof.





\section{Concluding remarks}\label{section:conclusion}

\subsection{Placements that correspond to the zero vector}

When $n=8k+1$ we can describe some placements that correspond to the zero vector in Nul $A$ as follows.  Select $k$ squares with coordinates $(x_1,y_1), \ldots, (x_k,y_k)$ such that these meet the following conditions: $0< x_i,y_i \leq 2k, \frac{y_i}{x_i} <1$, the values $x_1,y_1,\ldots,x_k,y_k$ are distinct and for any $i \neq j$ we have $\frac{y_i}{x_i} \neq \frac{y_j}{x_j}$. Next consider `folding' the board across the $x-$axis, then folding it across the $y-$axis, and finally folding it along the line $y=x$.  Consider the $8$ squares that `stack' on top of a selected square and place a queen in each of these.  Do this for each selected square.  We have now placed $8k$ queens.  Finally, place a lonely queen at $(0,0)$, for a total of $8k+1$ queens.  It is easy to check that this placement is in Nul $A$.  From a geometric perspective, the placement will look like $k$ octagons that are `nested' and centered at $(0,0)$.  The placement given in Figure \ref{fig:10on9by9withlonelycentered} is such a placement derived from this scheme where $n=9$ and the set of initial selected squares consists of one square with coordinates $(4,1)$; it is one of six possible placements for $n=9$ via this scheme.

When $n=5$ there are no placements that meet the conditions of Case 2 since a queen in $\QQ \setminus \QQ'$ forces the existence of four other queens and together with the lonely queen the placement would have at least $6>n$ queens.  When $n=8k+5$ we have not found any placements that correspond to the zero vector in Nul $A$.

\subsection{Other cases to consider}

When $n=4k+3$, we also believe that $m_3(n) \geq n+1.$  At this point of our investigations, we have not been able to establish this lower bound.  However, there are some sub-cases for which we can establish this improved lower bound: we sketch the proof here by following the polynomial method template provided in Case 1 of the proof of Theorem \ref{theorem:1mod4}.  By considering a good placement ${\cal Q}$ of size at most $4k+3$ on $\ZB{n}$ with at most $2k$ defined lines of slope $\pm 1$, the polynomial of degree $8k+4$

\begin{equation*}\label{eq:factorfor4k+3}
  g(x,y) = \prod_{j=1}^{2k+2} (x-\alpha_j)(y-\beta_j)\prod_{j=1}^{2k}(x-y-\gamma_j)(x+y-\delta_j)
\end{equation*}

can be constructed so as to vanish on all squares of the chessboard yet has a non-zero coefficient on the monomial  $x^{4k+2}y^{4k+2}$ by the Binomial Theorem.  Thus, by Theorem \ref{theorem:CN} we reach a contradiction.  This approach will settle all cases where $q' \geq 2$ and some others.  However, one particularly challenging sub-case that we cannot resolve is when $q'=0$ and the number of lines defined by ${\cal Q}$ is maximized (i.e. there are $2k+1$ lines defined by ${\cal Q}$ in each of the four slopes).  In this sub-case, one might consider the following polynomial  

\begin{equation*}\label{eq:factorfor4k+3-B}
  h(x,y) = \prod_{j=1}^{2k+1} (x-\alpha_j)(y-\beta_j)(x-y-\gamma_j)(x+y-\delta_j).
\end{equation*}

A leading monomial suitable for an application of the Combinatorial Nullstellensatz is again $x^{4k+2}y^{4k+2}$.  However, this monomial has the same coefficient as the same monomial in\\ $x^{2k+1}y^{2k+1}(x^2-y^2)^{2k+1},$ which is zero.  Thus, a polynomial method approach seems problematic.  At the same time, the combinatorial approach that comes towards the end of Case 2 in the proof of Theorem \ref{theorem:1mod4} seems problematic for the following reason.  In such a placement, there will be up to four queens that are neither lonely nor in a defining line of each possible slope.  Gaining information about the coordinates of these queens appears difficult and it could be that these particular queens are located on the perimeter of $U$.

Turning to a different set of cases: there are small cases of $n$ even where it has been established that $m_3(n)=n+1.$  From Table \ref{table:1}, we see that these values are $n=8, 14, 16, 20, 22, 24$.  There is no discernible pattern to us.

Recently, Di Stefano, Klav\u{z}ar, Krishnakumar, Tuite and Yero \cite{distefanoklavzarkrishnakumartuiteyero} have extended Gardner's problem to graph theory.  This, too, looks like an interesting line of inquiry.

\section{Acknowledgments}

The second author wishes to thank Fulbright Austria and the Algebra and Number Theory Group at the University of Graz for the support of his work and their warm hospitality.  This work commenced when he was a visitor there in the spring of 2022.

The authors wish to thank Greg Warrington of U. Vermont for reading an early draft and suggesting the idea of rotating the placement at the end of Case 2 in the proof of Theorem \ref{theorem:1mod4}.

\end{document}